\documentclass{article}
\usepackage{graphicx} 
\usepackage[a4paper, total={7in, 9in}]{geometry}
\usepackage{xcolor}
\usepackage[sort]{cite}
\usepackage{amssymb}
\usepackage{hyperref}
\usepackage{rotating}
\usepackage{amsmath}
\usepackage{tcolorbox}
\newtheorem{lemma}{Lemma}

\newtheorem{definition}{Definition}

\newtheorem{remark}{Remark}
\usepackage{booktabs} 
\begin{document}

\title{A novel statistical workflow for nonstationary modelling of successive Fr\'{e}chet extremes}
\date{}
\author{Grace Burtenshaw, Joe Lane, Meagan Carney}

\maketitle

\begin{abstract}

Accurate estimation of the frequency and magnitude of successive extreme events in energy demand is critical for strategic resource planning. Traditional approaches based on extreme value theory (EVT) are typically limited to modelling isolated extreme events and struggle to capture the dynamics of temporally clustered extremes, such as those driven by prolonged extreme weather events. These limitations are exacerbated by the scarcity of historical data and computational costs of longrun simulations leading to high uncertainty in return level estimates for successive extremes. Here, we introduce a novel statistical framework leveraging recent theoretical advances in successive extreme value modelling in dynamical systems. Under reasonable assumptions of the time series data (e.g. the data follow a fat-tailed Fr\'{e}chet distribution), our tool allows for significantly more robust estimates of returns and magnitudes of successive extreme events compared to standard likelihood methods. We illustrate our statistical workflow on scenarios of forecasted gas supply levels from 2025 to 2050. Common measures of statistical accuracy are provided as benchmarks for comparison. 
\end{abstract}

\section{Introduction and Background}

\subsection{Successive extremes in gas demand and the energy market in Australia}

Accurately estimating extremes in gas demand is crucial for safeguarding the financial interests of various industries and protecting the well-being of communities and individuals dependent on the energy market. Fluctuations in energy demand are patched through the use of gas generation. As such, uncertainties in energy usage can lead to unexpected extremes in gas demand. For example, there was a recent gas supply crisis occurring in Victoria attributed to unexpected extreme weather causing cold-spells (reduced solar energy to the grid) and a wind drought (reduced wind energy to the grid) throughout the 2024 winter \cite{wood2024victoria}. Given the growing instability and rapid change in Australia's energy systems, it is paramount to prioritise the development of methods that adequately quantify the likelihood of high impact, low-probability future gas demand events \cite{zuur2020resilient, elliott2018role}. This is critical for planning the decarbonisation of Australia’s National Electricity Market (NEM) which incorporates the interconnected electricity system across the five east-coast states of Queensland, New South Wales, Victoria, South Australia and Tasmania. The term NEM is commonly used interchangeably to represent both the physical electricity system as well as the market systems that drive the operation of that system. This study focuses on informing models of the former. 

Unlike western counterparts, Australia has a much higher reliance on variable renewable energy (VRE) influenced by lack of access to natural less-variable alternatives such as hydroelectric power, ambient heat and geothermal systems \cite{energygovau2024electricity, eurostat2024renewable}. Uncertainties in future weather variability, electricity demand, ongoing performance of conventional supply infrastructure (e.g. coal generators; transmission lines), and the stability and potential for future renewable electricity generation, support the need to investigate expected demand of gas as an alternative when renewable energy sources are insufficient.

The resulting combined variability, accompanied by the increasing likelihood of extreme weather events in recent years, contributes an additional complexity in modelling extremes in gas demand. These drivers force expected gas demand quantities to fluctuate, causing nonstationarity in the probability distributions of extremes in gas demand. Fortunately, this issue of nonstationarity can be overcome using more advanced statistical techniques from Extreme Value Theory (EVT) which allow the extremal model to change according to varying influences (or covariates). These techniques are well-known in the applied energy sector and have been used in recent literature to more accurately model extremes occurring in short bursts (or single peaks) \cite{chan2017extreme}; however, calls are still being made to incorporate more extremal statistics to improve the accuracy of models in energy systems \cite{natenergy2020, mccollum2020}. We refer the reader to \cite{coles2001, chan2017extreme} for a nice summary on current statistical techniques for modelling single occurrence extremes in time series.

Quantifying the impact of extreme events within the energy market has become a significant topic of discussion in contemporary literature. Understanding market trends and accurately forecasting events can be key to obtaining and maintaining a competitive edge when trading, developing market strategies, or designing policies. The authors in \cite{chan2017extreme} argued that by blending vector autoregression modelling (VAR) with standard extreme value methodologies, they were better able to capture volatility dynamics influenced by a variety of variables that impact the electrical market. Namely, seasonality, the tendency of the data to revert to an equilibrium level (often the marginal cost), negative pricing, and volatility clustering. 

Energy supply issues are compounded when extreme demand occurs over consecutive days or events such as in response to consecutive extreme weather over a region, e.g. heat-waves and cold-spells encouraging air-conditioning and heat usage  \cite{energycouncil, esci}.
A challenge often encountered when estimating extreme events is the scarcity of data on rare phenomena which increases uncertainty in the extremal model. This challenge is amplified for successive extremes, where the likelihood of observing runs of 2 or more extremes in a short time is significantly lower. Whilst improved forecasting results in energy demand have been achieved through the use of more classical methods from extreme value theory, these methods are limited to modelling extremes occurring in short bursts (or single peaks), rather than those occurring consecutively. Furthermore, the impact of successive extreme events in the context of gas-powered generation (GPG) demand remains widely unexplored. With the discussion of the impact of successive extreme climate events gaining momentum, given the inherent sensitivity of energy markets to weather, it is natural that research efforts extend their focus to understanding and modelling consecutive extremes in the energy space. This is particularly important for developing robust risk management strategies, improving market resilience, and ensuring the stability of energy systems in the face of increasing climate variability.

For the purposes of this investigation, we define successive extremes as the observation of consecutive high values occurring in a short run, for example, high levels of gas demand occurring over $k$ consecutive days. In the sections to follow, we develop and validate a statistical workflow that can be used to model returns and magnitudes of successive extremes in gas demand under varying drivers of non-stationarity. We illustrate the statistical robustness of our workflow through a numerical case-study using simulated data supplied by the UQ Energy Transition Modelling (ETM) team. These data represent a range of varying scenarios driving future gas demand over the years 2025 to 2050.

\subsection{Recent advancements in theoretical models of successive extremes}
Recent advancements in the field of dynamical systems have provided insights into key factors for more accurate modelling of successive extremes. To explain the impact of these results, we begin by introducing some background on the relationship between extreme value theory for dynamical systems and its role in deriving extreme value statistics for gas demand.

In the dynamical framework, one assumes that the state of a complex system (e.g. the state of the energy grid) can be represented by some coupled system of equations (e.g. the weather, population, government initiatives, etc.) within which one can measure an observable of interest (e.g. gas demand) that changes according to changes in the system. To gain insights into the expected behaviour of an observable it is common practice to measure and record at (often equally-space) time intervals, where $X_i$ corresponds to the value of the observable at time interval $i$. The sequence $(X_n)$ that is generated from this process is called a (nonlinear) time-series and can be treated as a stochastic process wherein theoretical insights in probability and statistics can be invoked to make predictions on the future or long-term behaviour of values of $(X_n)$.

\begin{remark} In the real world, we often do not know the underlying system of equations; however, we are still able to generate a nonlinear time-series of observations by manually measuring the observable of interest through field tools.
\end{remark}

To study the behaviour of extremes of $(X_n)$, we consider the sequence of maxima (or minima) defined by,
\[
M_n = \max\{X_1,\dots,X_n\}.
\]
Beginning with work of \cite{leadbetter1983} and \cite{chernick1981} in the 1970s, probability distributions for $(M_n)$ have been mathematically proven to exist for a wide-range of sequences $(X_n)$ under assumptions of weak dependence (satisfying mixing and recurrence conditions $D(u_n)$ and $D^q(u_n)$) and stationarity (the distribution of $(X_n)$ does not change according to some covariate or time). This limiting probability distribution of $(M_n)$ is called an Extreme Value Law. The formal definition of the Extreme Value Law is provided below.

\begin{definition}[Extreme Value Law] Let $(u_n)$ be a sequence of constants defined by the requirement that $$\lim_{n\rightarrow\infty} n\mathbb{P}(X_1>u_n)=\tau,$$ and suppose the sequence $(X_n)$ is stationary and satisfies $D(u_n)$ and $D^k(u_n)$ then,
\[
\lim_{n\rightarrow\infty}\mathbb{P}(M_n\le u_n) = e^{-\theta\tau}
\]
\end{definition}
where the parameter $\theta\le 1$ is called the extremal index and $1/\theta$ (roughly) measures the expected cluster size of extremes.

For numerical applications of data modelling, the sequence $(u_n)$ is taken as a sequence of increasing thresholds under which a frequency can be calculated, such that $u_n= ya_n+b_n$ where $(a_n)$ and $(b_n)$ are normalising sequences with $(a_n )>0~\forall n$. For some fixed $n$, the sequence $(a_n)$ and $(b_n)$ take fixed values in $\mathbb{R}$ which can be understood as some scale $\sigma$ and location $\mu$ parameter. Under this framework, we can rewrite the Extreme Value Law in the following way:
$$P(a_n (M_n-b_n )\le y)\rightarrow G(y)$$
where $G(y)$ is the limiting extremal distribution of the suitably normalised maxima and takes the form of one of three types of distributions. 

A Type 1 (Gumbel) distribution has an upper exponentially decaying tail; Type 2 (Fr\'{e}chet) corresponds to an upper polynomial decaying tail; and Type 3 (Weibull) corresponds to a bounded tail. Work from \cite{fisher1928} and \cite{gnedenko1943} guarantee that all three extreme value distributions may be combined into the following three parameter distribution function,
\begin{definition} [Generalized Extreme Value Distribution (GEV)]
$$P(M_n\le z) \rightarrow G(z) = \exp\bigg(-\bigg[1+\xi \bigg(\frac{z-\mu}{\sigma}\bigg)\bigg]^{1/\xi}\bigg)$$
where $\mu$ is the location parameter, $\sigma$ is the scale parameter, and $\xi$ is the shape parameter which determines the tail-behaviour or type of distribution \cite{coles2001}. 
\end{definition}

The weak dependence criteria $D(u_n)$ holds for sequences $(X_n)$ with sufficiently fast (polynomial) decay of correlations - this assumption is almost always guaranteed in real-world systems - while $D^q(u_n)$ holds for sequences $(X_n)$ with extremes that do not occur too frequently in large succession. On the other hand, the assumption of stationarity is often violated in real-world systems. We discuss this problem and a well-accepted statistical solution at length for our purposes in the sections to follow.

Provided the data meet the assumptions of stationarity and weak dependence, the limiting distribution of $(M_n)$ is known, and creating a numerical model for $(M_n)$ is then reduced to applying an appropriate statistical fitting method that results in the most accurate numerical approximation of the parameters $\xi$, $\mu$, and $\sigma$. Since traditional statistical fitting methods require a frequency plot of the data to estimate parameters, the block maxima method is typically employed to generate such a frequency plot in which the data is blocked into pieces of length $m$ (long enough to guarantee convergence to the GEV), the maxima are calculated over each block, and this sequence of maxima is used to approximate the parameters of the corresponding GEV. 

Recent advancements in understanding and modelling successive extremes have centred around taking functions the sequence $(X_n)$, for example the moving minimum over some short $k$-length window and interpreting returns of extreme highs or lows of this function as indicative of $k$ successive high or low values of $(X_n)$. This can be seen by noting that if $Y_i$ is above some high threshold, then all the values in the $k$-length window must also be above threshold \cite{carney2023}. We now restate for our context the dynamical lemma which will be central to our proposed statistical workflow for accurate modelling of successive events (and days) of extreme gas demand across Australia.

\begin{lemma}[Lemma 4.1 \cite{carney2023}]\label{scaling}
Suppose $M_n = \max\{X_1,\dots,X_n\}$ follow a GEV, $G(\xi_1,\mu_1,\sigma_1)$ of Frech\'{e}t type. Define $Y_j = \min\{X_j,\dots,X_{j+k-1}\}$ for $j = 1,\dots n$ to be the minimum value over a time window of length $k$ and $B_n = \max\{Y_1,\dots Y_n\}$, so that if $B_i\ge z$ for $z\in\mathbb{R}$ then every value over the window of length $k$, ensures $X_i\ge z,\dots,X_{i+k-1}\ge z$. Suppose $(X_n)$ and $(Y_n)$ are stationary, weakly-dependent sequences, satisfying $D(u_n)$ and $D^q(u_n)$, with extremal index $\theta_1$ and $\theta_k$, respectively. Then the maximal process $(B_n)$ follows a GEV, $G(\xi_k, \mu_k,\sigma_k)$ of Frech\'{e}t type with,
\[
\xi_k = \xi_1, \quad \mu_k = g_T(k)\mu_1\bigg(\frac{\theta_k}{\theta_1}\bigg)^{\xi_1}, \quad \sigma_k = g_T(k)\sigma_1\bigg(\frac{\theta_k}{\theta_1}\bigg)^{\xi_1}
\]
for some unknown function $g_T(k)$ depending on the window size $k$ and the system $T$. 
\end{lemma}

From the scaling lemma we conclude that if we have an extreme value law preserved for sequences of maxima $(M_n)$ and $(B_n)$, then the scaling function $g_T(k)$ can be used to derive the parameter estimates for the distribution of $(B_n)$ from those of the corresponding distribution of $(M_n)$. The authors in \cite{KardkasemCarney} invoke the scaling lemma to aid in the choice of limiting distribution by shape stability for successive extremes in rainfall across Australia. We build on this concept by creating a workflow that (1) requires a fixed shape across all windows of length $k$ to aid in faster convergence of the limiting distribution for $(B_n)$ and (2) numerically approximates $g_T(k)$ to infer the scale and location of the limiting distribution for $(B_n)$. In general, our workflow will allow users to estimate the probability distribution of $k$-successive extremes using the probability distribution of a single extreme where much more data is available.

\section{Methodology}
\subsection{Scaling function $g_T(k)$ as a generalised linear model}\label{gkt}
The goal of this section is to rewrite the scaling function $g_T(k)$ as a generalised linear model (glm), so that it can be readily estimated using statistical fitting methods. A good approximation of $g_T(k)$ from the data allows us to infer the scale and location parameters of the extremal probability distribution for predictions on return times and magnitudes of any number of $k$-successive extremes along the sequence $(X_n)$.

From the scaling lemma \ref{scaling}, we note that the derived relationships allow us to conclude the location (similarly, the scale) for any number of $k$-successive extremes relates to the single extreme case by
\begin{align*}
\mu_k &= g_T(k)\mu_1\bigg(\frac{\theta_k}{\theta_1}\bigg)^{\xi_1}.
\end{align*}
In the cases to follow it is important to note that, for some reasonable choices of $k$ upon which we fit our $g_T(k)$ model, values of the extremal index for $(X_n)$, $\theta_1$, and $(Y_n)$, $\theta_k$, can be computed through well-known robust, numerical methods such as Ferro-Segers \cite{ferro2003} and values of $\mu_k$ and $\mu_1$ can be computed via likelihood methods. There is a prediction horizon on $k$, after which we can no longer reliably compute these parameters due to lack of historical data and must invoke the scaling lemma relationships to obtain results. We explore this in detail and compare results to current statistical methods in section \ref{casestudy}.

\noindent\textit{Case 1. $g_T(k)$ follows a linear or polynomial relationship to $k$.}\label{gkt}
If $g_T(k) = a_1+a_2k+a_3k^2+\dots+a_mk^m$, it is straight forward to rearrange the scaling relationship outline in lemma \ref{scaling},
\begin{align*}
a_1+a_2k+a_3k^2+\dots+a_mk^m = \frac{\mu_k}{\mu_1\bigg(\theta_k / \theta_1\bigg)^{\xi_1}}
\end{align*}
and fit a glm using x-axis predictors $k,k^2,\dots k^m$.

\noindent\textit{Case 2. $g_T(k)$ follows a power-law relationship to $k$.}
We have $g_T(k) = a k^{\beta}$ so that,
\begin{align*}
\log(\mu_k) &= \log(g_T(k))+\log(\mu_1)+\xi(\log(\theta_k)-\log(\theta_1))\\
\log(\mu_k) &= \log(ak^{\beta})+\log(\mu_1)+\xi(\log(\theta_k)-\log(\theta_1))\\
\log(\mu_k) &= \beta\log(k)+\log(a)+\log(\mu_1)+\xi(\log(\theta_k)-\log(\theta_1)).
\end{align*}
Fitting a glm to the paired data $(\log(k), \log(\mu_k))$ results in $\beta$ approximated by the slope of the glm and $\log(a)+\log(\mu_1)+\xi(\log(\theta_k)-\log(\theta_1))$ approximated by the intercept. Furthermore, since $\mu_1$, $\theta_k$, and $\theta_1$ are computable, one can explicitly work out the value of $a$ to complete the formula representation for $g_T(k)$.

\noindent\textit{Case 3. $g_T(k)$ follows an exponential relationship to $k$.}
Here we have $g_T(k) = ab^{k-1}$ so that
\begin{align*}
\log(\mu_k) &= \log(g_T(k))+\log(\mu_1)+\xi(\log(\theta_k)-\log(\theta_1))\\
\log(\mu_k) &= \log(ab^{k-1})+\log(\mu_1)+\xi(\log(\theta_k)-\log(\theta_1))\\
\log(\mu_k) &= (k-1)\log(b)+\log(a)+\log(\mu_1)+\xi(\log(\theta_k)-\log(\theta_1)).
\end{align*}
Fitting a glm to the pair data $(k-1,\log(\mu_k))$ results in $\log(b)$ approximated by the slope of the glm and $\log(a)$ worked out in a similar fashion to \textit{Case 2}.

\begin{remark}
The relationship in Case 3. has been shown to hold for Fr\'{e}chet distributions of $(X_n)$ representing historical rainfall values over Australia and for general classes of chaotic dynamical systems with uniform expansion. This relationship cannot be assumed for general $(X_n)$; however, it is interesting that we find a similar result holds here for $(X_n)$ representing gas demand.
\end{remark}

\subsection{Nonstationary environments and adaptations for $g_T(k)$}\label{nongkt}
Nonstationarity in the probability distribution for the sequence $(X_n)$ is very common occurrence in real-world settings where drivers force the probability of an event to change. Some well-known examples of this are: the El Ni\~{n}o - La Ni\~{n}a cycles forcing hot/dry or cold/wet conditions over Eastern Australia; greater climate variability over time with rising sea levels and sea surface temperatures; and increased urbanisation over time causing increases in overall electricity demand. This nonstationarity causes the parameters of the extreme value distribution (most often the location and scale parameters) to change according to some driving force. Let $t$ be a covariate that changes the parameters of the distribution for $(X_n(t))$. The following examples are worked out for the location parameter while the scale is reserved for validating the approximation of $g_T(k)$ in this setting.

\noindent\textit{Example 1. The location parameter of the distribution for a single maxima changes linearly with $t$.}
In this setting, the location parameter of the distribution for a single maxima is represented by,
\[
\mu_1(t) = \mu_{1,0}+\mu_{1,1}t
\]
One can then use standard likelihood fitting methods to find the values $\mu_{1,0}$ and $\mu_{1,1}$ so that the extremal distribution statistically significantly fits the data. We refer the reader to \cite{coles2001} for a very nice overview of how this approximation can be achieved. Relating the nonstationary relationship to scaling lemma \ref{scaling} gives
\begin{align*}
\mu_k(t) &= g_T(k)\mu_1(t)\bigg(\frac{\theta_k(u_t)}{\theta_1(u_t)}\bigg)^{\xi_1}\\
\mu_{k,0}+t\mu_{k,1} &= g_T(k)(\mu_{1,0}+\mu_{1,1}t)\bigg(\frac{\theta_k(u_t)}{\theta_1(u_t)}\bigg)^{\xi_1}\\
\mu_{k,0}+t\mu_{k,1} &= g_T(k)(\mu_{1,0})\bigg(\frac{\theta_k(u_t)}{\theta_1(u_t)}\bigg)^{\xi_1} + g_T(k)(\mu_{1,1}t)\bigg(\frac{\theta_k(u_t)}{\theta_1(u_t)}\bigg)^{\xi_1}.
\end{align*}
The above holds true for any $t\in\mathbb{R}$ so we can fix $t = 0$ to obtain,
\[
\mu_{k,0} = g_T(k)(\mu_{1,0})\bigg(\frac{\theta_k(u_t)}{\theta_1(u_t)}\bigg)^{\xi_1}.
\]
Methods from section \ref{gkt} can now be easily implemented to approximate $g_T(k)$.

\noindent\textit{Example 2. The location parameter of the distribution for a single maxima changes exponentially with $t$.}
The approach in this example follows a very similar process to the above and is included to aid the reader in understanding how these methods would readily expand to include a vast number of examples in nonstationarity. Here, we assume the location parameter of the distribution for a single maxima is represented by,
\[
\mu_1(t) = \exp(\mu_{1,0}+\mu_{1,1}t).
\]
Invoking the scaling lemma \ref{scaling} gives
\[
\exp(\mu_{k,0}) = g_T(k)(\exp(\mu_{1,0}+\mu_{1,1}t))\bigg(\frac{\theta_k(u_t)}{\theta_1(u_t)}\bigg)^{\xi_1},
\]
which holds true for any $t\in\mathbb{R}$, so fixing $t = 0$ gives the relationship
\[
\exp(\mu_{k,0}) = g_T(k)(\exp(\mu_{1,0}))\bigg(\frac{\theta_k(u_t)}{\theta_1(u_t)}\bigg)^{\xi_1}
\]
and methods from section \ref{gkt} apply almost immediately, with the straightforward change that the approximated $g_T(k)$ predicts the exponential of $\mu_{k,0}$ and hence would need to be converted via the natural log to complete the analysis.

\subsection{Approximating the Extremal Index $\theta$ in a nonstationary environment}\label{extremal}
The process described above to approximate $g_T(k)$ in a nonstationary environment relies on an approximation of the extremal index $\theta_1(u_t)$ (for single maxima) and $\theta_k(u_t)$ (for $k$-successive maxima) as functions of the covariate $t$ describing the driving force of nonstationarity in the parameters of the generalised extreme value distribution. To understand how we address this issue, we first recall the definition of the extremal index as the parameter in the extreme value law that describes the expected cluster size of extremes in a short window \cite{leadbetter1983}. 

It has been shown for sequences that satisfy an extreme value law modelling \textit{magnitudes} of returns of extremes, a compound Poisson distribution can be used to model the \textit{number} of returns of extremes in the sequence \cite{lucarini2016}. Under the assumption of a compound Poisson limiting distribution for returns, the authors in \cite{hsing1988} have shown that $1/\theta$ can be interpreted as the expected cluster size of exceedances. That is, when one exceedance occurs, one would expect that $1/\theta$ exceedances on average occur in a short window. There is extensive literature on the interpretation and approximation of the extremal index. We refer the reader to \cite{lucarini2016} for a nice summary of the current methods used for numerical approximation of the extremal index. 

Given its proven numerical robustness, we have chosen to approximate the extremal index via the Ferro-Segers estimate \cite{ferro2003}. The estimate is well-accepted in the extremal community and, for replicability of results, is performed through a function that is built into many extremal platforms (e.g. MATLAB, R, Python). The Ferro-Segers estimate assumes a limiting compound-Poisson distribution for returns in our daily precipitation time-series and approximates the extremal index through a maximum likelihood calculation of the mean intensity of the compound-Poisson distribution defined by,
\[
\theta = 2\frac{\sum_{j = 1}^{N}(T_j)^2}{N\sum_{j=1}^N(T_j)(T_j-1)}
\]
where $T_j$ is the wait time associated to the exceedance $Z_j>u$ for $j = 1,\dots N$ over some high (fixed) threshold, often taken as the 95\% or 99\% quantile of the data $(X_i)$ for $i = 1,\dots n$.

We adapt this definition of the Ferro-Segers approximation for the extremal index to the nonstationary setting by defining the threshold as $u(t)$ where $t$ is the covariate describing the driving force of extremes. The choice of function $u(t)$ is informed by historical values of $(X_n(t))$ over $t$ where $u(t)$ is defined as some function of the 95\% or 99\% quantiles of the data over $t$. A simple and quick approach is to define $u(t)$ as the moving quantile of $(X_n(t))$ over some fixed window along $t$. We refer the reader to section \ref{casestudy} for case-study examples in this type of numerical approximation.

\subsection{Complete proposed statistical workflow}
From the derivations described above, we propose the following novel statistical workflow to model successive extremes of the Fr\'{e}chet type.
\begin{tcolorbox}[width=\linewidth, colback=blue!5!white, colframe=blue!60!black, title=Proposed Statistical Workflow]
\noindent Part A. Known statistical workflow for single extremes \cite{coles2001}.
\begin{itemize}
\item[] Step A1. Evaluate sources of nonstationarity in $(X_n)$ and investigate the autocorrelation of the sequence $(X_n)$ to inform the block maxima method.
\item[] Step A2. After appropriate choices of covariates and nonstationary models of the parameters $\mu_1$ and $\sigma_1$, fit a generalised extreme value (GEV) distribution on the sequence of maxima $(M_n)$ calculated over \textit{long enough} blocks of $(X_n)$ using standard likelihood methods.
\item[] Step A3. Use the shape parameter approximation of the GEV on $(M_n)$ to determine whether the extremes are of the Fr\'{e}chet type (e.g. $\xi_1>0$).
\end{itemize}
\tcblower
\noindent Part B. Proposed statistical workflow for successive extremes of Fr\'{e}chet type.
\\\\
\noindent For values of $k$ with sufficient available historical data for numerical approximation:
\begin{itemize}
\item[] Step B1. Calculate the moving minimum $(Y_n)$ of $(X_n)$ over window size of length $k$ defined in lemma \ref{scaling}.
\item[] Step B2. Compute the maxima $(B_n)$ of $(Y_n)$ over blocks the same length as Part A.
\item[] Step B3. With a \textbf{fixed} shape $\xi_k = \xi_1$, fit the GEV with the same nonstationary models of the parameters $\mu_k$ and $\sigma_k$ as in Part A. **Even with a fixed shape, likelihood fitting will eventually fail for large enough $k$, where convergence rates to the GEV are reduced as observed by the authors in \cite{carney2023}.**
\item[] Step B4. Compute the nonstationary extremal index $\theta_1$ and $\theta_k$ for all $k$ using the modified Ferro-Segers definition described in section \ref{extremal}.
\item[] Step B5. Fix covariate values (e.g. $t = 0$) in the fitted nonstationary parameters of the GEV for each $k$ and use all appropriate data pairs (e.g. $(k-1,~ \log(\mu_{0,k})$), or similar depending on the structure of $g_T(k)$) to estimate the function $g_T(k)$ in the nonstationary setting as described in section \ref{gkt} and section \ref{nongkt}
\end{itemize}
\noindent Part C. Use the relationship from lemma \ref{scaling} and the numerically approximated $g_T(k)$ to model the GEV (and hence, returns and magnitudes) for any number of $k$ successive extremes of $(X_n)$.
\end{tcolorbox}

\subsection{Alternative approach for extremes with a nonstationary shape parameter}
In the proposed workflow we assume, perhaps presumptuously, that the shape parameter $\xi_1$ for $(M_n)$ does not suffer from nonstationarity. This assumption suggests that the tail decay behaviour is not changed via some underlying force. Fortunately, this process can be applied in the same way, provided the distribution of $(M_n)$ does not change in type. This can be observed by noting that since max-stability extends for any number of $k$ successive extremes from lemma \ref{scaling} we may take
\[
\xi_k(t) = \xi_1(t),
\]
as the 'fixed' shape in Step B3 without issue.

It should still be noted that there is evidence that some real-world data may suffer from this phenomenon (see, for example \cite{forster2025} illustrating extreme rainfall tail decay changes according to geographic location). In this case, a unified approach for any extremal type would be ideal and certainly a future research direction of interest for the authors.

\section{Numerical case-study on simulated gas demand}\label{casestudy}
\subsection{Gas demand model}
The Energy Transition Modelling (ETM) team of the UQ Gas and  Energy Transition Research Centre (UQ-GET) provided the data for this project \cite{uqget}. The ETM initiative is building capability for assessing uncertainties in the modelling of long-term energy system transition pathways for Australia. The purpose of this case study was to contribute to an active UQ-GET project that is evaluating risks to future NEM resilience, associated with uncertainty over the ability of the east-coast gas supply system to meet an increasingly volatile demand for GPG in the future. Here our work in returning probabilities for consecutive, rare events aims to lessen and quantify uncertainty in future demand for gas-powered generation.  

The data was generated by UQ-GET, simulating the evolution of the NEM through 2050 by capturing a range of uncertainties over future weather variability and different future trajectories for infrastructure deployment. Specifically, underpinning the two scenarios detailed in this paper were varied schedules for infrastructure modification. In scenario p2-002, there were no proposed modifications to the infrastructure. Comparatively, the simulation for scenario p2-016 accounts for 1-3 year delays in the timing of new pumped-hydro projects across various Australian states. In developing and applying extreme value analysis to those results, this application aims to demonstrate techniques that would improve the quantification of planning risk. The modelling environment is comprised of multiple layers of stochastic partial differential equations whose parameters are influenced by historical weather data. These nuances introduce randomness in the system, requiring a statistical approach to understanding expected returns in gas demand.

Each scenario in the model returns the forecasted duration, total energy and average energy demand (defined as the average amount of energy demand in terajoules (TJ) per day in a single gas demand event) for gas events occurring from 2025 to 2050. A gas event is defined as a point in time where gas is supplied to the market. We focus on extremes in average daily energy demand from gas generation. Due to the nature of gas demand being supplementary to the energy market, it is common to observe large bouts of time with zero gas demand. For illustrative purposes and to avoid very long block requirements, we remove all gas events where zero demand occurs and perform our statistical workflow to predict the likelihood of $k$ successive extreme average energy demand events. We illustrate our technique in the following section on two simulated gas demand scenarios.

\subsection{Building the reference distributions for $(M_n)$ and $(B_n)$}
For appropriate scaling choices, the exact limiting distribution is the generalised extreme value distribution (GEV), $G(\xi,\mu,\sigma)$, with shape $\xi$, location $\mu$ and scale $\sigma$ parameters. In practice, the \textit{block maximum approach} is commonly employed which involves dividing the time-series into blocks of a fixed length $m$ and performing some likelihood fitting of the parameters of the GEV using the sequence of block maxima $(M_{n_m})$. Blocks are chosen long enough so that correlations between $X_{n}$ and $X_{n+m}$ have sufficiently decayed and to guarantee $(M_{n_m}$) is in the limiting regime of the GEV. Quality of fits on the sequence of maxima $(M_{n_m})$ for the approximated $(\hat{\xi},~\hat{\mu},~\hat{\sigma})$ are often initially investigated through quantile plots. Applying the standard method of maximum likelihood estimation, we find very poor fits for the block maxima of the data under the assumption that the distribution of gas demand, $(X_n)$ is stationary (see Appendix Figure \ref{fig:non}). Prompted to investigate possible nonstationarity in the data, we find time serves as a driving covariate in the parameters of the GEV with relationships summarised in Table \ref{tab:relationships}.

\begin{table}[h!]
\centering
\begin{tabular}{|c|c|c|c|}
\hline
\multicolumn{4}{|c|}{Nonstationary Model and $g_T(k)$ Type by Scenario} \\
\hline
\multicolumn{2}{|c|}{Scenario p-002}  & \multicolumn{2}{|c|}{Scenario p-016}\\
\hline
Nonstationary Model & $g_T(k)$ Type & Nonstationary Model & $g_T(k)$ Type \\
\hline
$\mu_k(t) = \mu_{k,0}+\mu_{k,1}t$ & & $\mu_k(t) = \mu_{k,0}+\mu_{k,1}t$ & \\
$\sigma_k(t) = \exp({\sigma_{k,0}+\sigma_{k,1}t})$ & $ab^{k-1}$ & $\sigma_k(t) = \sigma_{k,0}+\sigma_{k,1}t$ & $ab^{k-1}$ \\
$\xi_k = \xi_k$ & & $\xi_k = \xi_k$ & \\
\hline
\multicolumn{2}{|c|}{Likelihood ratio $p$-value: $<10^{-5}$}  & \multicolumn{2}{|c|}{Likelihood ratio $p$-value: $<10^{-5}$} \\
\hline
\end{tabular}
\caption{Data-derived relationships to be used in the proposed statistical workflow.}\label{tab:relationships}
\end{table}

Following standard methods in the field \cite{coles2001}, we fit the nonstationary GEV on the data $(M_{n_m})$ via maximum likelihood estimation. Estimates are provided in Table \ref{tab:nonstationaryGEV} where $k = 1$ corresponds to the distributional parameters of $(M_{n_m})$. These new fits result in markedly better quantile fits. Kolmolgorov-Smirnov and Anderson-Darling tests were also performed and indicate statistically significant fits at the 95\% confidence level. Further details on nonstationary fitting and illustrations can be found in the Appendix \ref{append}.

Next, we implement Step B1. - Step B3. of the proposed statistical workflow for $k = 1,\dots, k_f$, assuming a fixed shape for every window size $k$, estimated from $k = 1$, $(M_{n_m})$ where the maximum amount of data on extremes is available. The fitted parameters are then used to build the reference distributions for $(B_{n_m})$, noting that the distribution for $(B_{n_m})$ is only valid up to some fixed window size $k\le k_f$ for which there is sufficiently available successive extreme data to numerically fit the distribution. This relationship will be used to derive the necessary function $g_T(k)$ that will allow us to infer distributional parameters for $(B_{n_m})$ for any window size, including those beyond estimate horizon $k_f$.

\subsection{Numerical estimation of the scaling function $g_T(k)$ and parameter inference}
\begin{figure}[h!]
\centering
\includegraphics[width = 0.5\textwidth]{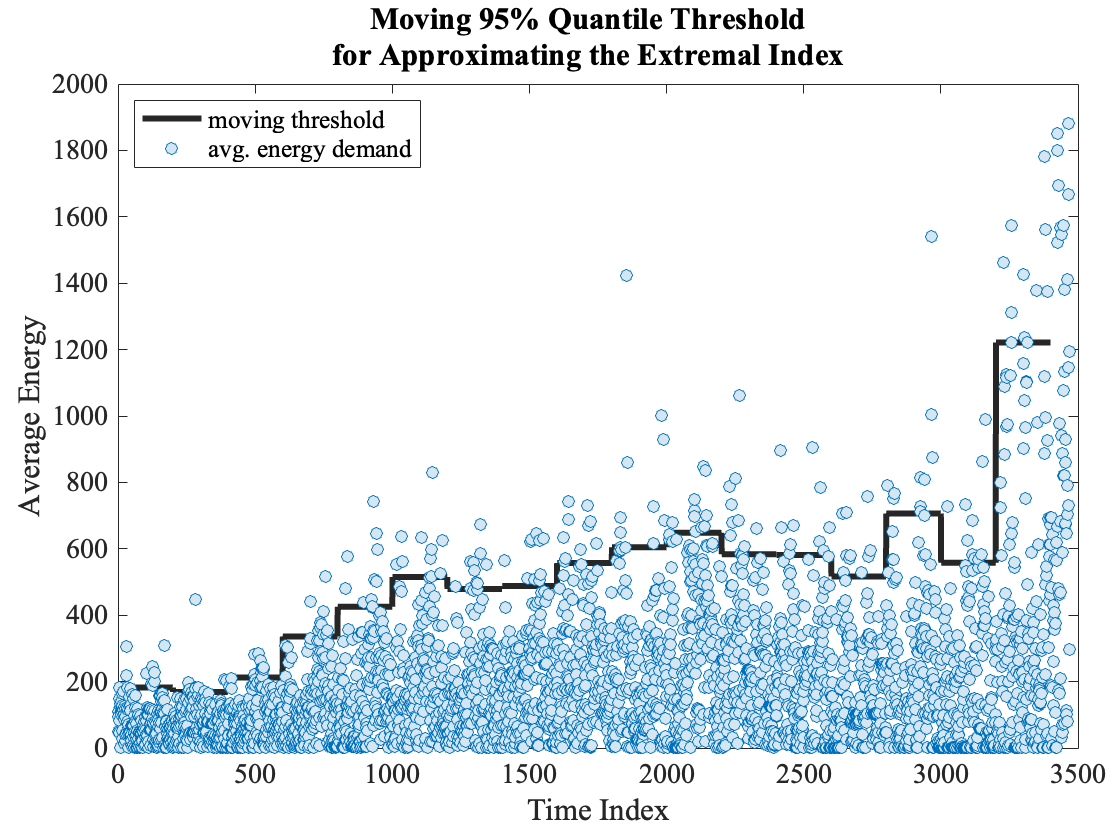}
\caption{Illustration of nonstationary thresholding for extremal index approximation. Here, we use the moving yearly 95\% quantile as the threshold for both scenarios.}\label{fig:theta}
\end{figure}

Addressing step B4. of the proposed statistical workflow, we numerically estimate $\theta_1$ and $\theta_k$ in the nonstationary setting through the adapted Ferro-Segers approximation where thresholds are chosen to be the 95\% quantile for each year of extreme gas demand. Figure \ref{fig:theta} illustrates the chosen nonstationary thresholding for extremal index approximation. Estimates of the extremal index over increasing window size $k$ are provided in Table \ref{tab:ei}.

From the maximum likelihood estimated  values of $\mu$, $\sigma$ and $\xi$ for window sizes  $k = 1,\dots, k_f$ resulting from steps B1. - B3. of the proposed statistical workflow, we numerically investigate the relationship between the value of $\mu$, for some fixed $t$, and the window size $k$ following the described workflow and examples from Section \ref{nongkt}. We find a linear relationship on the pair $(k-1,\log(\mu))$ indicating $g_T(k) = ab^{k-1}$ according to derivations in Section \ref{gkt}, \textit{Case 3}. A bit of algebra is required to find the exact values $a$ and $b$ in the model of $g_T(k)$ where we note that, for the generalised linear model fit on these pairs we return a slope $m$ and intercept $\beta$ for the model $y = mx+\beta$ where $y = \log(\mu_k)$, $x=k-1$, $m = \log(b)$ and $\beta = \log(a)+\log(\mu_1)+\xi(\log(\theta_k)-\log(\theta_1))$, hence, 
$$b = \exp(m)$$ 
and 
$$a = \exp(\beta-\log(\mu_1)-\xi(\log(\theta_k)-\log(\theta_1))).$$
It is straight-forward to observe that for $\mu_1$ with a linear relationship to $t$, we need only take $t = 0$ and fit the data on the pair $(k-1, \log(\mu_{k,0}))$ where $\mu_{k,0}$ is the intercept of the nonstationary model. Our method can be extrapolated to any choice of nonstationary model, as discussed in Section \ref{nongkt}, with minor modifications. We illustrate this for our case-study in the validation phase on the scale parameter $\sigma$ in the following paragraphs.

Once the parameters of $g_T(k)$ are numerically fit using generalised linear model methods on the location parameter, we use the resulting $g_T(k)$ to predict values of $\sigma$ for the same fixed $t$. Exactly how this is done is determined by the nonstationary model of the scale parameter. We illustrate this by example for each of the scenarios in our case-study.

\begin{figure}
\centering
\includegraphics[scale=0.90, trim = {20, 375, 0 0} ]{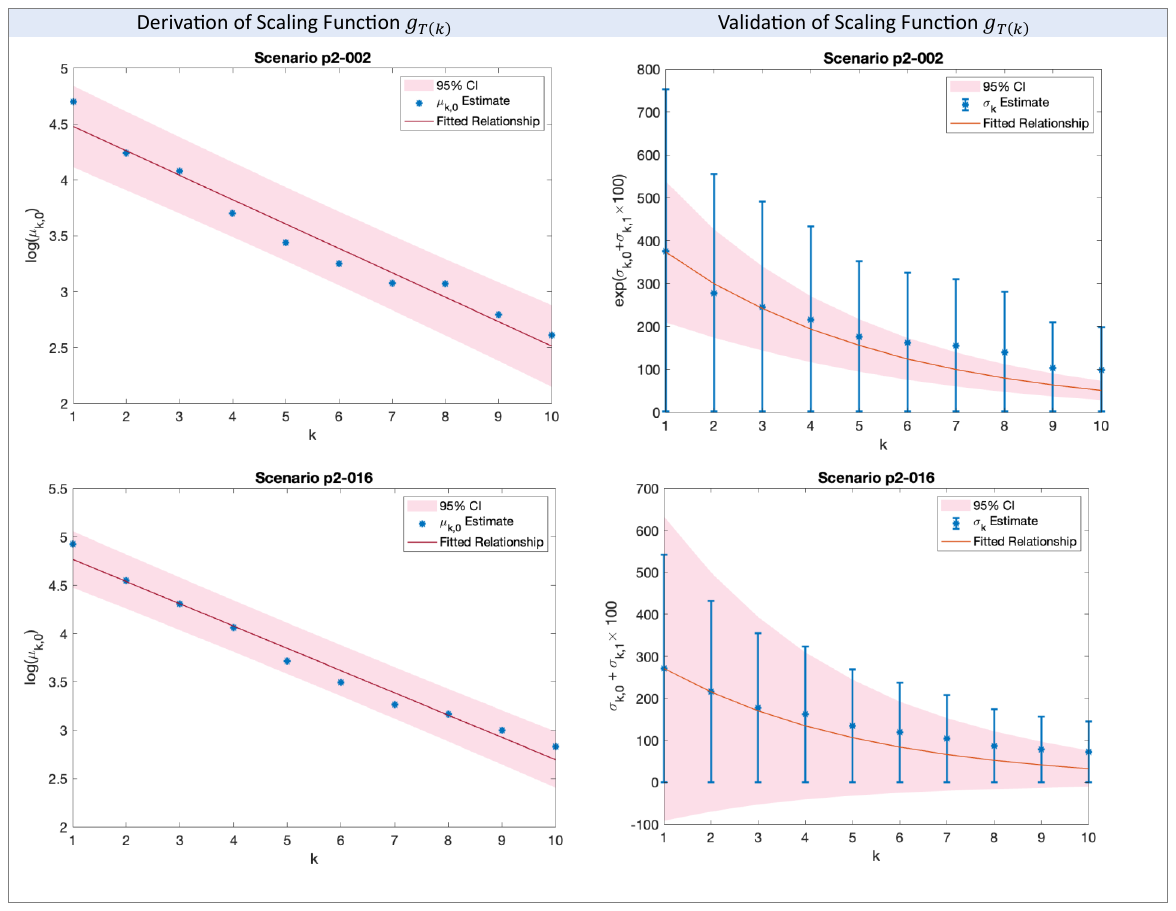}
\caption{Derivation and validation plots of the scaling function $g_T(k)$. 95\% confidence intervals on the fit of $g_T(k)$ and the maximum likelihood estimate of the parameter are indicated in red shading and blue error bars, respectively.}\label{fig:gkt}
\end{figure}

\noindent\textit{Scenario p-002:} The scale parameter has the relationship $\sigma_k(t) = \exp(\sigma_{k,0}+\sigma_{k,1}t)$ so that for fixed $t = 0$ we have, 
\begin{align*}
\exp(\sigma_{k,0}) &= g_T(k)\exp(\sigma_{0,0})\bigg(\frac{\theta_k}{\theta_1}\bigg)^{\xi_k}.
\end{align*}
\noindent\textit{Scenario p-016:} The scale parameter has the relationship $\sigma_k(t) = \sigma_{k,0}+\sigma_{k,1}t$ so that for fixed $t=0$ we have,
\begin{align*}
\sigma_{k,0} &= g_T(k)\sigma_{0,0}\bigg(\frac{\theta_k}{\theta_1}\bigg)^{\xi_k}.
\end{align*}
Derivation (on $\mu$) and validation (on $\sigma$) complete the investigation that makes up step B5. in the proposed workflow. Figure \ref{fig:gkt} illustrates the results of this derivation and validation against the true values of $\mu_k$ and $\sigma_k$, respectively, for fixed $t=0$ (e.g. $\mu_{k,0}$ and $\sigma_{k,0}$). Estimates from generalised linear modelling of $g_T(k)$ are reported in Table \ref{tab:gkt}.

\subsection{Improved performance over traditional fitting methods}\label{sec:improved}
\begin{figure}[h!]
\centering
\includegraphics[width = 0.5\textwidth]{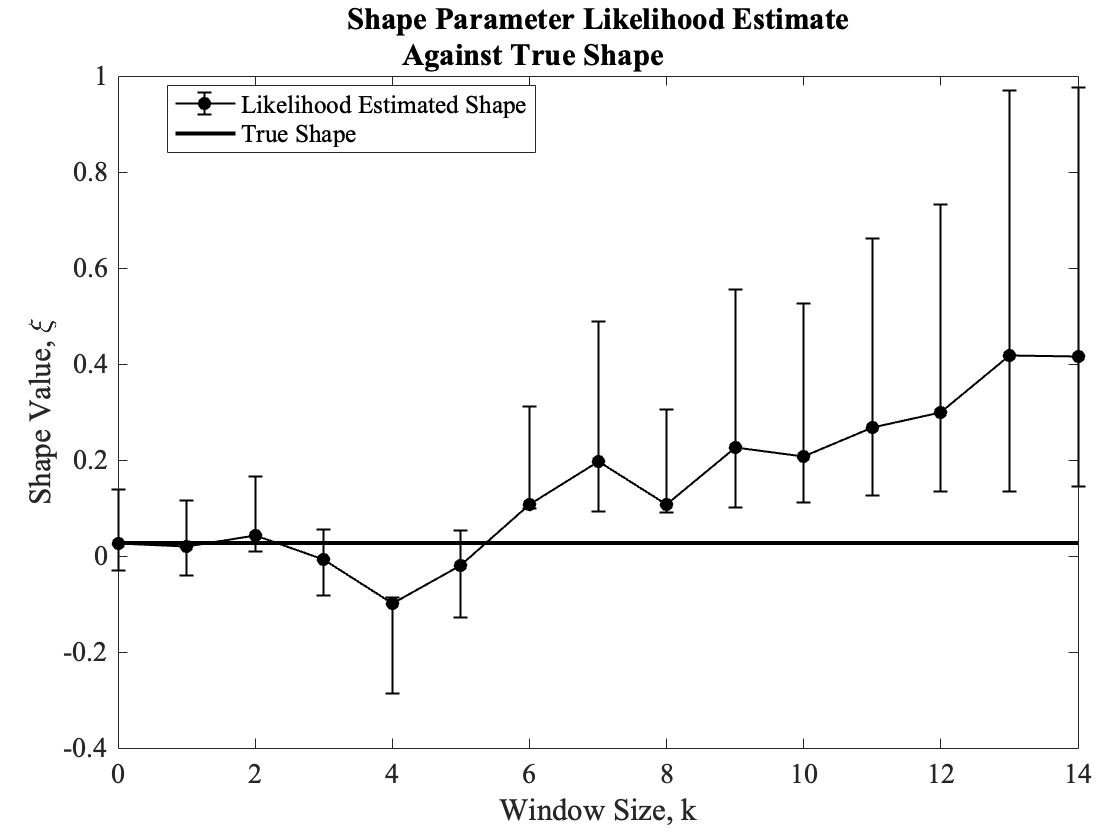}
\caption{Illustration of shape parameter estimates using likelihood estimates that rely on the quantity of available rare event data. As the number of successive extremes (window size) increases, and hence becomes more rare, the likelihood estimate of the shape parameter deviates from its true value.}\label{shape}
\end{figure}

\noindent\textit{The numerical importance of fixing the shape parameter in step B3.} For increasing window sizes (e.g. increasing numbers of successive extremes) we observe a poor estimate of the shape parameter due to lack of available tail data to fit the distribution of $B_n$. Furthermore, poor estimates of the shape parameter increase error in the number estimates of the location and scale parameters. This phenomenon is illustrated in Figure \ref{shape}.

\begin{figure}[h!]
\centering
\includegraphics[width = \textwidth, trim = {0, 350, 0 0}]{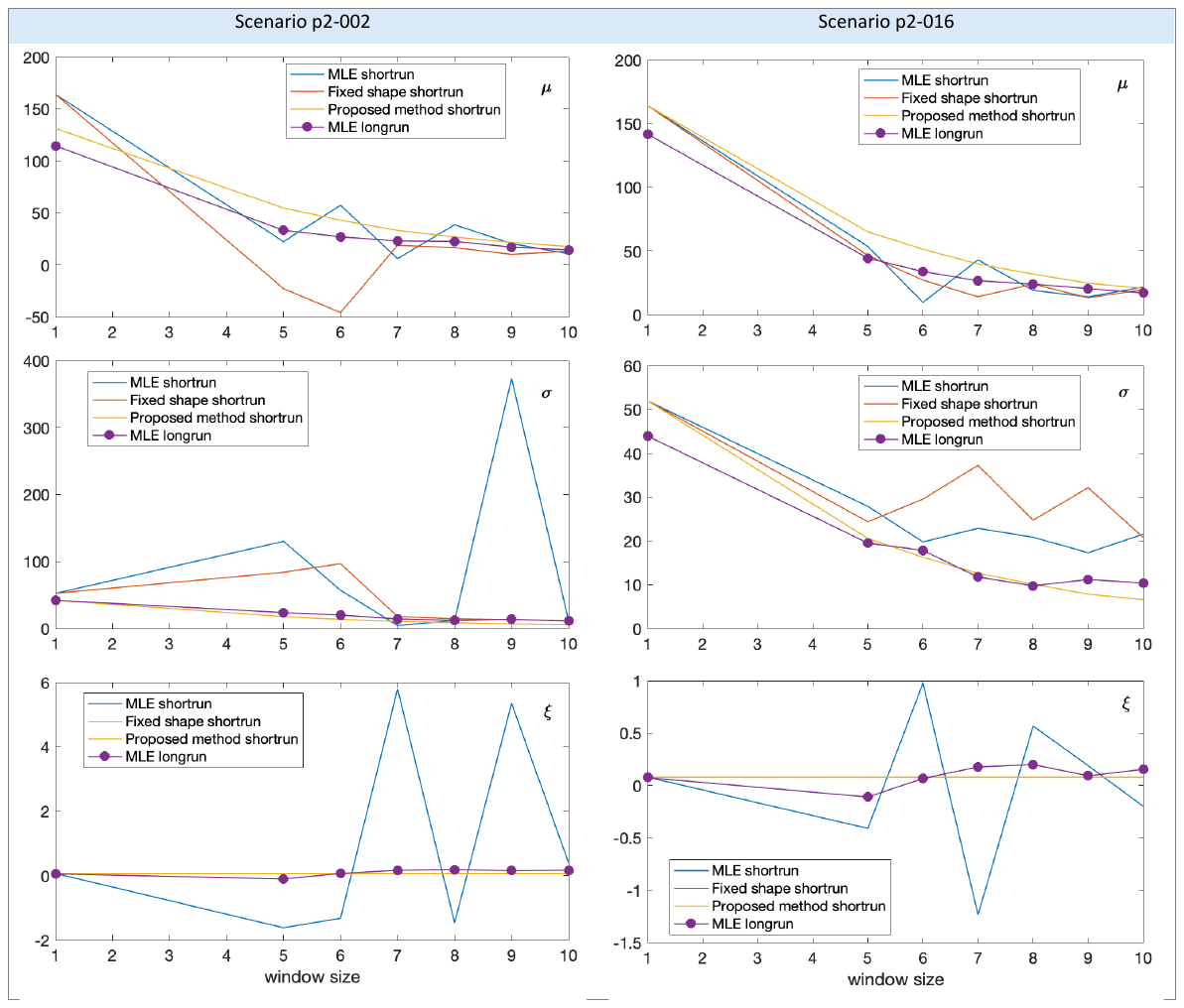}
\caption{Parameter estimates of the extremal distribution for fixed time $t = 1$ using each of the methods (e.g. the current modelling framework based on maximum likelihood estimation, an assumed fixed shape proposed here as an improvement, and the full proposed method proposed here as the ideal modelling framework). Purple marked lines indicate the long-run parameter estimates representing the true distributional parameters.}\label{fig:params}
\end{figure}

\noindent\textit{The numerical improvement of parameter inference over maximum likelihood estimation for $k>k_f$.}
With a fixed shape, we estimate $g_T(k)$ on values of location and scale parameters coming from lower window sizes $k\le k_f$ and then use this information to infer the values of the parameters for window sizes beyond the numerical estimation horizon $k>k_f$. To numerically test the accuracy of our method against traditional likelihood techniques, we first estimate the \textit{true} distribution using a long-run of 3,347 historical datapoints of non-zero gas demand. Then, we perform estimation of the extremal distribution parameters on a short-run of 300 historical datapoints of non-zero gas demand using,
\begin{itemize}
\item[(1)] Our full proposed method to approximate each of the parameters of distribution (ideal modelling framework proposed here)
\item[(2)] The maximum likelihood method with an assumed fixed shape parameter (improved modelling framework proposed here)
\item[(3)] The maximum likelihood method only (current modelling framework for successive extremes)
\end{itemize}
We compare these estimates to the parameters from the long-run distribution. Figure \ref{fig:params} shows the estimates and true values of the parameters for a fixed time $t = 1$, noting that this plot will be the same (but shifted) distribution for any choice of $t$. Scenarios p2-002 and p2-016 indicate significantly better performance for our full proposed method (1) against the fixed shape assumed MLE method (2) and the MLE method (3). Importantly, we note that confidence intervals cannot be reliably computed for each window size $k$ (e.g. number of successive extremes) with such low sample size (300 points).

\begin{figure}[h!]
\centering
\includegraphics[width = \textwidth]{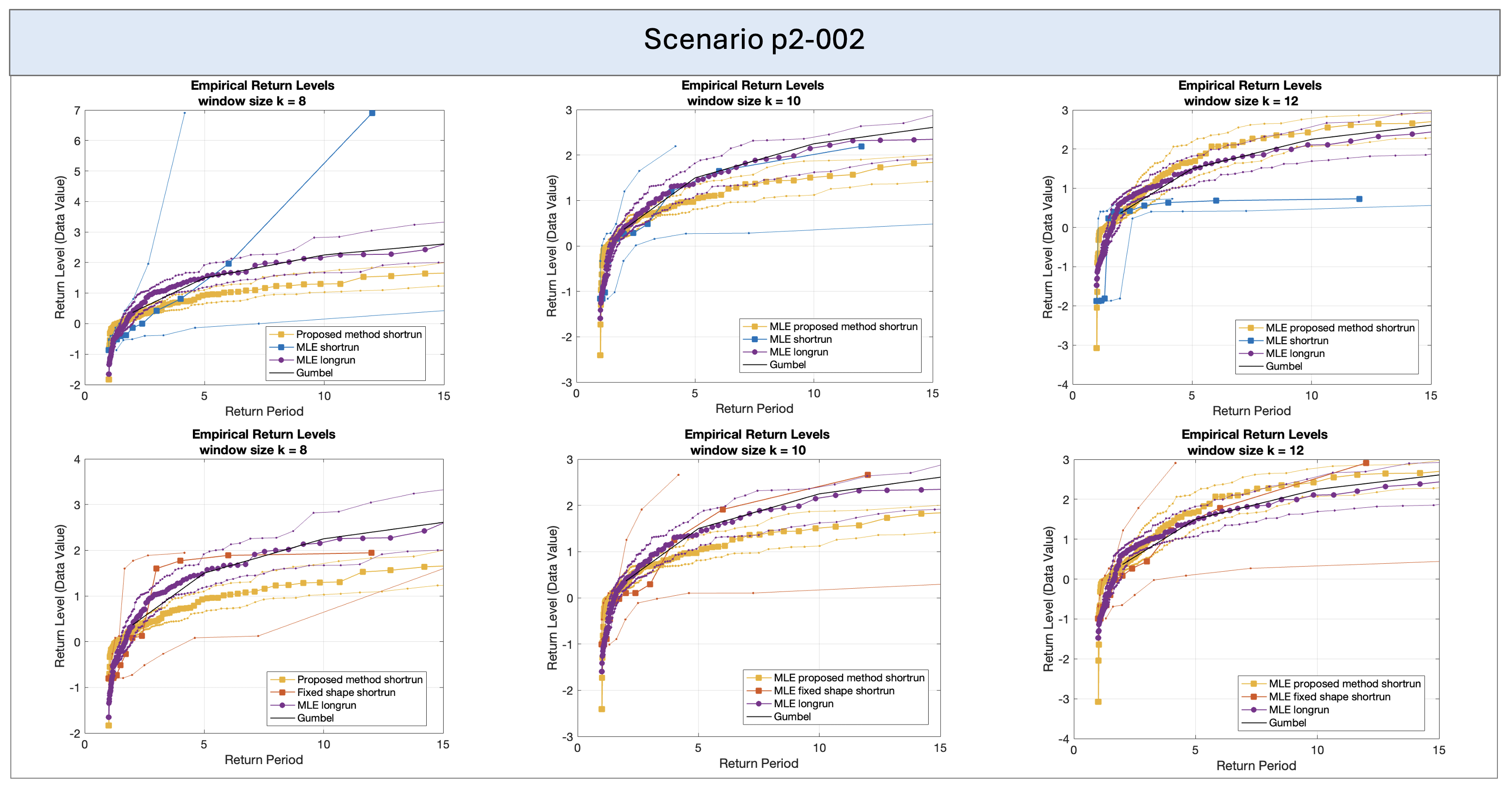}
\includegraphics[width = \textwidth]{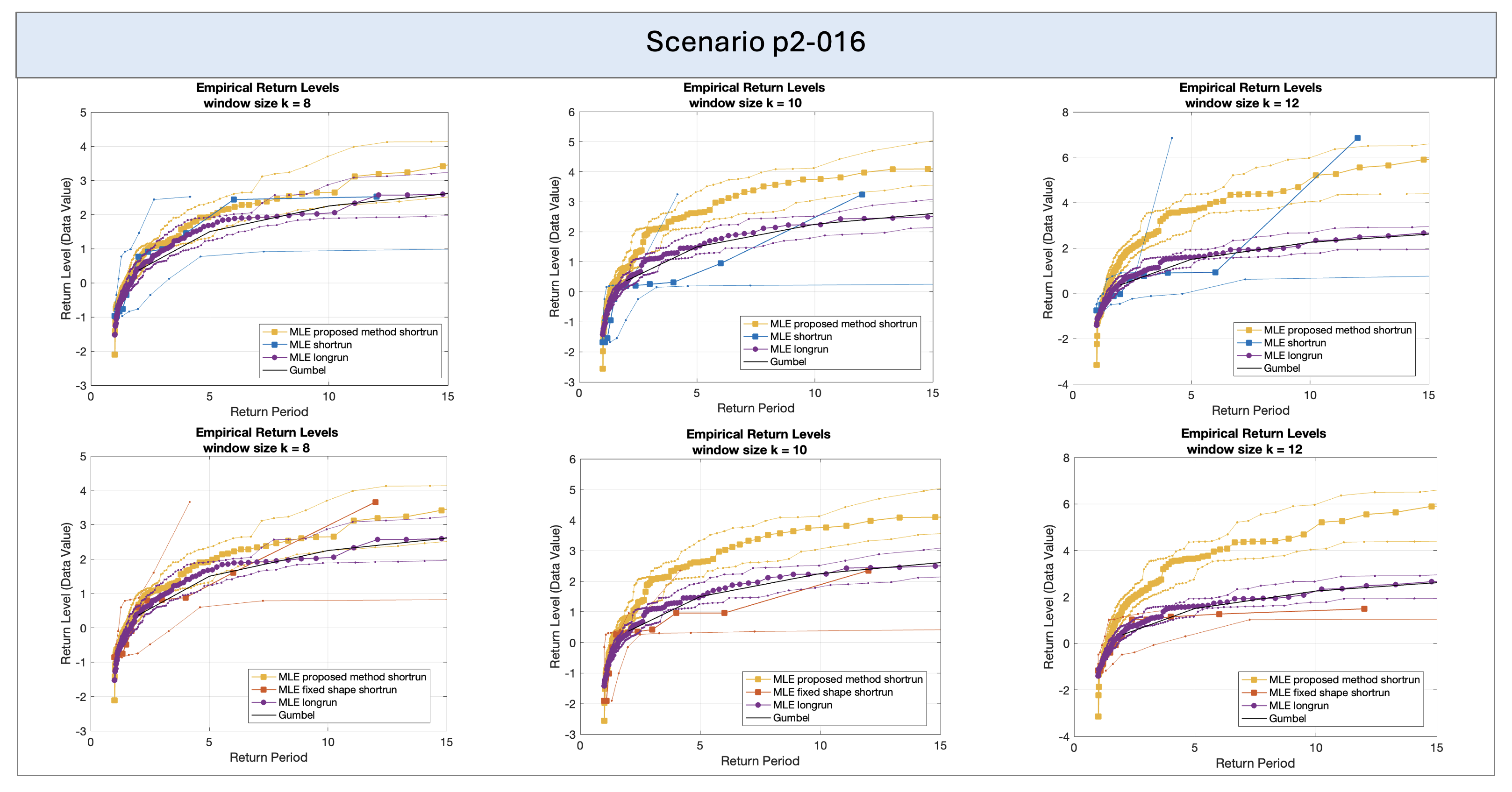}
\caption{Return level - return time plots from the proposed method compared to standard likelihood methods that make up the current modelling framework in successive extremes. This is illustrated using both methods on each simulated scenario for window sizes $k = 8$, $k = 10$ and $k = 12$, modelling the likelihoods of observing eight, ten and twelve successive over-threshold values (extremes) of gas demand, respectively. Transformations to the standard Gumbel distribution were performed given the nonstationary nature of the distribution.}\label{fig:rlrt}
\end{figure}

Next, we generate and evaluate the performance of the return level-return time plots generated from each of the methods (1), (2), and (3) and compare these results against the return level-return time plot of the long-run distribution. The issue of nonstationarity in illustrating the fit of each method is overcome by normalising the data to the standard Gumbel distribution as described by the authors in \cite{coles2001} using the parameter estimates from each of the three methods. We find two important numerical results: first, methods (2) and (3) relying on the likelihood evaluation at each window size $k$ fail to generate estimates for return times earlier than the long-run and our full proposed method (1) and second, our full proposed method (1) outperforms methods (2) and (3) for increasing return times and increasing window size. We illustrate these results in Figure \ref{fig:rlrt}. Where possible, 95\% confidence intervals were generated using the built-in bootstrap approach available for empirical cdf approximation in MATLAB. Missing confidence lines for increasing return periods indicate numerically unreliable approximations due to low sample size. Scenario p2-002 reflects longer periods of reliability for our method; however, this is certainly situational as observed in p2-016 where error increases significantly over return period for longer windows despite better overall performance against traditional likelihood fitting. We investigated the reason for this increase in error and found that it is due to a small difference in the approximation of the scale parameter which is pronounced for larger window sizes (since the scale is significantly smaller) and is further emphasised in the required standard Gumbel normalisation for the return time-return level plot.

\subsection{Estimating return probabilities of successive extremes in gas demand}
\begin{figure}[h!]
\centering
\includegraphics[width = 0.5\textwidth]{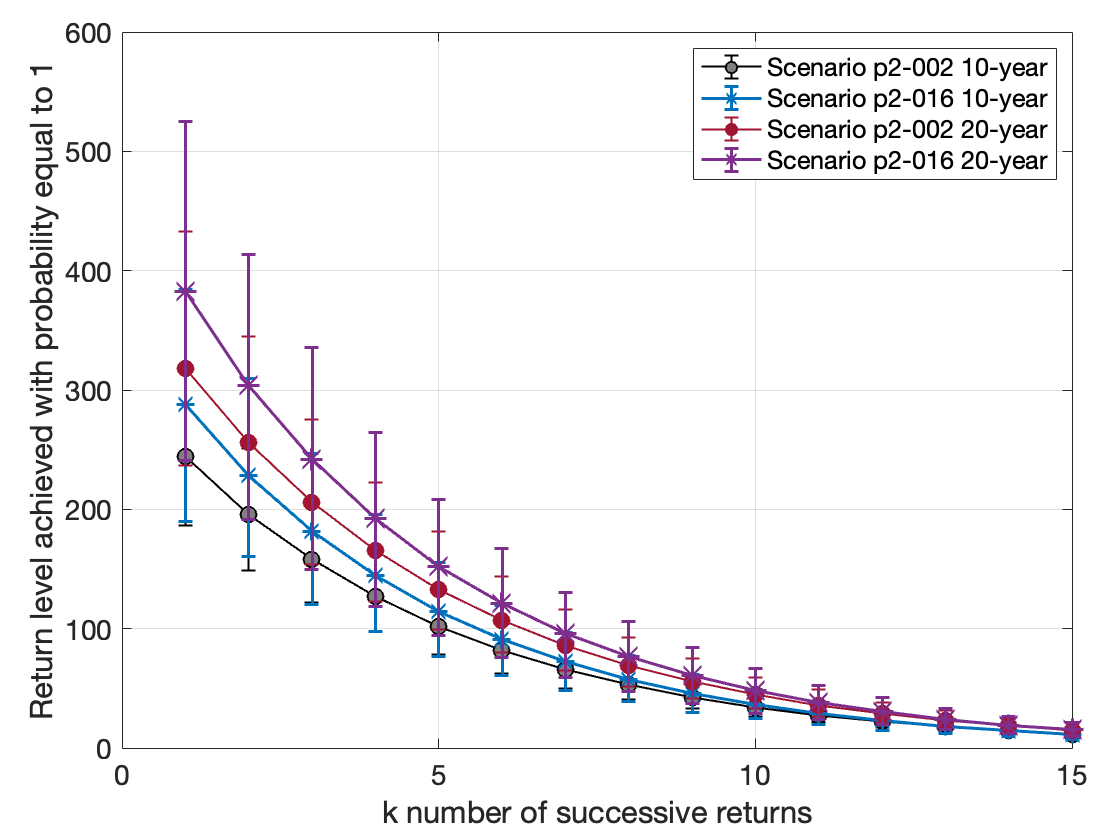}
\caption{Return levels and corresponding 95\% confidence intervals found through permutations. Reported return levels are the values of successive maximum gas demand we expect exceed in 10- and 20-years with probability equal to 1. Plots are illustrated for each scenario along increasing successive windows $k$.}\label{fig:results}
\end{figure}

\begin{figure}[h!]
\centering
\includegraphics[width = \textwidth]{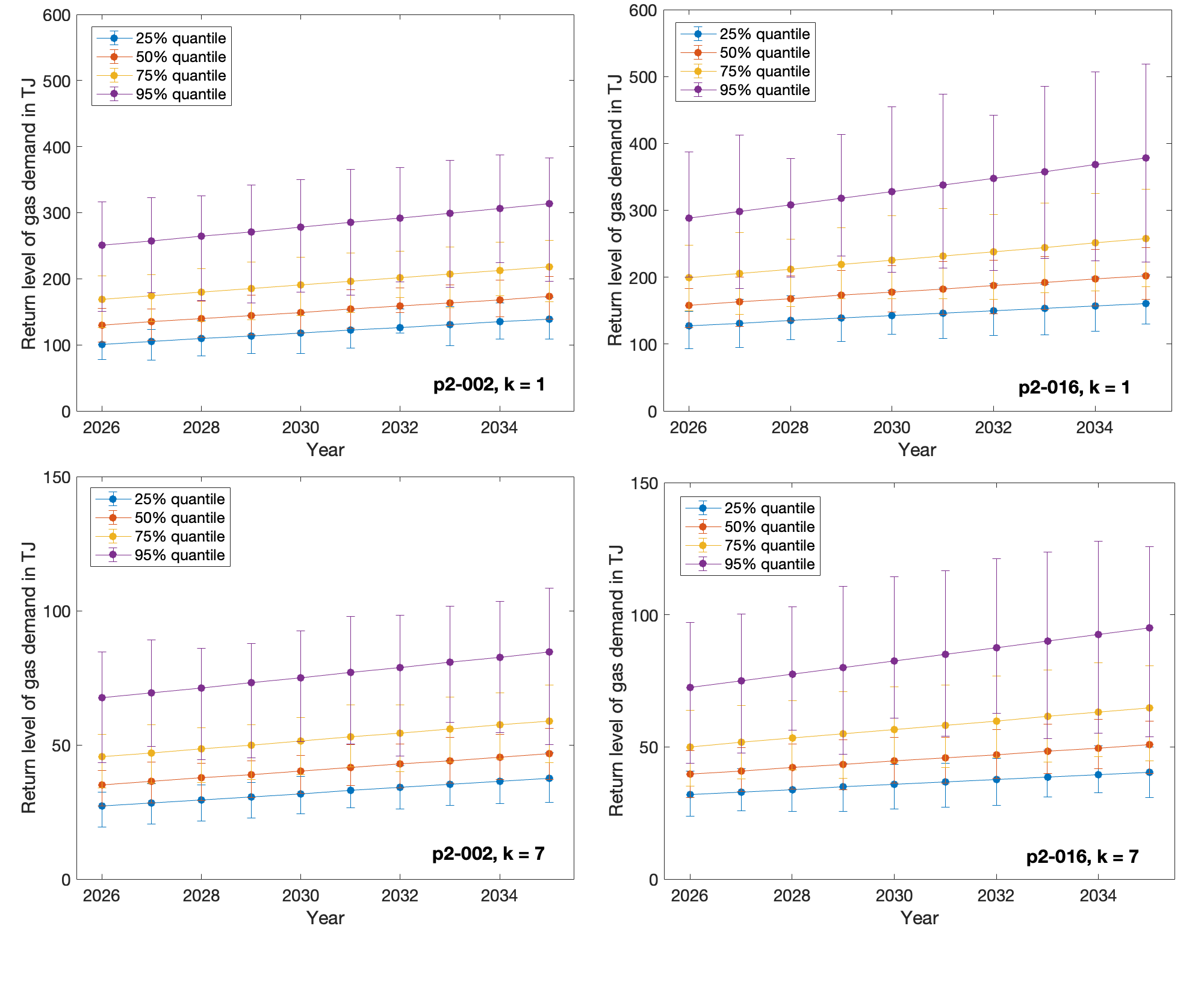}
\caption{Return levels of single occurrence extreme gas demand and $k = 7$ successive returns of extreme gas demand corresponding to 25\%, 50\%, 75\% and 95\% quantiles with confidence intervals for each scenario (p2-002 and p2-016) estimated from probability distributions generated from our proposed method for each year in the next 10 years (2026-2035).}\label{fig:quants}
\end{figure}

We complete our analysis by illustrating how our proposed method can be used to more accurately approximate returns of successive extremes in gas demand under a nonstationary driving force. Using our final model from the proposed method with parameters and scaling function $g_T(k)$ fit over all available data, we first estimate the $k$ successive return levels that we expect Australia will exceed with probability one. Using these values, we can specify the $k$ successive quantity of gas demand we expect will occur at least once in Australia over the next 10 and 20 years, respectively. These levels are provided in Figure \ref{fig:results} with 95\% confidence intervals for $k = 1,\dots,15$ successive extremes. Confidence intervals over the quantiles were created by drawing 10 samples from the 95\% confidence intervals for each of the estimated parameters $\{\mu_{k,0}~\mu_{k,1},~\sigma_{k,0},~\sigma_{k,1},~\xi\}$, computing the probability distributions for every permutation of the samples, and returning the minimum and maximum of the corresponding quantiles across all the probability distributions. Interestingly, the 10-year / 20-year return level estimates indicate that a single extreme gas demand event is likely to occur beyond a 243.88TJ / 318.05TJ threshold in scenario p2-002 (no change in the current policy) and beyond a 287.74TJ / 382.92TJ threshold for scenario p2-016 (1-3 year delays in hydropower infrastructure). Moreover, a one-week successive extreme gas demand event is likely to occur beyond a 65.81TJ / 85.83TJ threshold for scenario p2-002 and a 72.28TJ / 96.20TJ threshold for scenario p2-016. That is, gas demand is expected to occur beyond these thresholds for 7 successive events, or 460.67TJ / 600.81TJ in total for scenario p2-002 and 505.96TJ / 673.40TJ in total for scenario p2-016.  These results indicate that, in simulated scenarios the policy change marked in scenario p2-016 is likely to result in higher demand thresholds; however, their differences diminish for successive extremes indicating only a small difference of successive extreme event exceedances occurring beyond $k = 7$. Furthermore, the total quantity for $k = 7$ suggests that in these scenarios we are likely to observe a 7-day successive extreme gas demand event that requires 79\% (88\%) of the peak supply capacity of the Woleebee Creek Facility, the largest Queensland gas reserve, within the next 20 years in scenario p2-002 (p2-016) \cite{aemo_gas_bulletin_board}. Interestingly, we also note from Figure \ref{fig:quants} that increased uncertainty is observed for scenario p2-016, indicating that the proposed 1-3 year delays in hydropower infrastructure, according to the UQ-GET simulation may contribute significantly increased uncertainty in extreme gas demand (e.g. 95\% quantile).

We end by providing key quantile approximations for returns of successive gas demand projected to occur over the next 10 years (e.g. 2025 - 2035). Figure \ref{fig:quants} illustrates increased demand levels in time across all $k$ successive extreme gas demand events. That is, all things remaining the same, the data from the UQ-GET model indicate that Australia is likely to observe higher exceedances of successive gas demand events over time. This observation coincides with the latest interim report on gas demand published by the Australian Competition and Consumer Commission (ACCC), in which they use the related rolling 7-day average demand rather than the successive 7-day demand observable used in this investigation. Without intervention, the ACCC quarterly report indicates that supply may become unstable in 2025 due to increased overall GPG demand contributing to shortfalls in reserves; however, optimistically this report also highlights that residential, commercial and industrial demand for gas are expected to remain steady due to increase electrification, further highlighting the benefits of policy intervention to decrease blackout risk even in the context of increasing uncertainty. The ACCC report is completed with a call for methodologies needed to incorporate weather-related variation in projections of demand. This work, along with the weather-driven simulations provided by the team from UQ-GET takes a significant step in answering this call.

\section{Discussion}
We propose a novel statistical workflow to model successive Fr\'{e}chet extremes and illustrate that this workflow outperforms conventional extreme value modelling methods. The workflow is built on the premise of generating a new sequence as a function of the original sequence, the windowed moving minimum, that preserves the successive property of extremes by way of noting that if the minimum over a window is above some value, then every observation in the window much be above that value. Preserving the consecutive nature of the data comes at the cost of data scarcity since the likelihood of observing $k$ extremes successively decreases as the number of successions, $k$, increases. Such scarcity reduces our ability to fit the data using conventional extreme value modelling methods. As a first improvement, we emphasise the extension of the max-stability criterion recently proven by the authors in \cite{carney2023} and proposed the method of fixing the shape parameter for the extremal distribution of increasingly successive $k$ extremes of a sequence. Although seemingly simple, this improvement guarantees that we estimate the shape parameter, the slowest converging parameter of the generalised extreme value distribution, using the original sequence where data is least scarce. We further improve the fits of our distribution for increasing $k$ by numerically approximating a function, called $g_T(k)$ that is used to estimate the last two parameters (location and scale) of the generalised extreme value distribution. Unlike in the stationary dynamical setting of \cite{carney2023}, we must address adaptations to the nonstationary setting in a rigorous way to guarantee a working statistical methodology for real-world settings.

We illustrate the effectiveness of our proposed methodology on simulations of gas demand from UQ-GET. Using the resulting non-stationary generalised extreme value distributions for $k$ successive extreme gas demand events, we explored the flexibility and robustness of the methodology by applying the workflow to two infrastructure scenarios p2-002 and p2-016. Importantly, we find that the UQ-GET simulations indicate increased extreme gas demand and significantly uncertainty under the assumption of a 1-3 year delay in hydropower infrastructure development. These improvements to demand forecasting and uncertainty can provide stakeholders and policymakers with more effective tools to optimise facility and reserve capacity metrics, increase market security in the face of unprecedented demand, and inform frameworks and policy which drive investments in infrastructure resilience.

We emphasise that the models introduced here are informed by the UQ-GET model of simulated GPG demand. Consequently, one very important caveat to this analysis is that these data used to inform the models in this analysis do not take into consideration any government or public policy changes that may occur in future years. For example, the Clean Energy Council of Australia has proposed several alternatives to gas demand in the form of solar or alternative renewable energy replacements \cite{energygovau2024electricity}. As such, this review may serve as a risk evaluation to aid in quantity of renewable energy storage that may be needed in future infrastructure to compensate for projections of increased \textit{successive} energy demand events whose uncertainties of occurrence are a source of difficulty in achieving net zero goals.

Importantly, whilst we illustrate our proposed methodology on successive extremes in gas demand, it can be used in a much more general context to model successive extremes in nonlinear timeseries in both stationary and non-stationary environments, provided the data follow a heavy-tailed (Fr\'{e}chet) distribution.

\subsection{Limitations}
This study is subject to limitations. The most obvious is the requirement for the data to follow a Fr\'{e}chet distribution. Whilst this investigation does not cover the Weibull or Gumbel limiting regimes, it is certainly an active area of interest for the authors who are currently working on a unifying theory through the extensions of lemma similar to Lemma \ref{scaling}. 

Another limitation is the prediction horizon of $g_T(k)$ whose accuracy we expect will diminish with increasing $k$. It is a fine numerical balance between having enough estimates of the location and scale along $k$ to approximate the function $g_T(k)$ whist not inferring the estimate of location and scale for $k$ too large due to increasing error. Apriori error approximations on the estimates using the confidence intervals of $g_T(k)$ made be made through bootstrapping approaches or via the delta method, as employed here.

\clearpage

\begin{table}[h!]
\centering
\begin{tabular}{|c|c|c|c|c|c|c|c|c|c|c|c|c|c|c|c|c|}
\hline
\multicolumn{11}{|c|}{Fitted Values (Maximum Likelihood Estimated)} \\
\hline
$k$ & \multicolumn{5}{|c|}{Scenario p-002}  & \multicolumn{5}{|c|}{Scenario p-016}\\
\hline
& $\mu_{k,0}$ & $\mu_{k,1}$ & $\sigma_{k,0}$ & $\sigma_{k,1}$ & $\xi_k$ & $\mu_{k,0}$ & $\mu_{k,1}$ & $\sigma_{k,0}$ & $\sigma_{k,1}$ & $\xi_k$\\
\hline
1 & 110.11 & 4.49 & 3.72 & 0.02 & 0.06 & 137.58 & 4.24 & 42.25 & 1.72 & 0.08 \\
\hline
2 & 69.39 & 3.19 & 3.44 & 0.02 & 0.06 & 94.36 & 2.75 & 27.62 & 1.42 & 0.08 \\
\hline
3 & 58.94 & 2.43 & 3.13 & 0.02 & 0.06 & 74.32 & 2.11 & 14.26 & 1.23 & 0.08\\
\hline
4 & 40.44 & 2.08 & 3.11 & 0.02 & 0.06 & 58.05 & 1.66 & 13.13 & 1.12 & 0.08 \\
\hline
5 & 31.17 & 1.60 & 3.12 & 0.02 & 0.06 & 41.13 & 1.30 & 16.87 & 0.88 & 0.08 \\
\hline
6 & 25.83 & 1.23 & 2.97 & 0.02 & 0.06 & 32.90 & 1.05 & 17.06 & 0.76 & 0.08 \\
\hline
7 & 21.66 & 1.01 & 2.66 & 0.02 & 0.06 & 26.22 & 0.87 & 11.08 & 0.70 & 0.08 \\
\hline
8 & 21.52 & 0.70 & 2.46 & 0.02 & 0.06 & 23.80 & 0.64 & 8.92 & 0.59 & 0.08 \\
\hline
9 & 16.36 & 0.63 & 2.59 & 0.02 & 0.06 & 20.08 & 0.59 & 10.64 & 0.51 & 0.08 \\
\hline
10 & 13.62 & 0.52 & 2.37 & 0.02 & 0.06 & 16.92 & 0.46 & 9.66 & 0.47 & 0.08 \\
\hline 
\multicolumn{11}{|c|}{$1.96\times\sigma/N$ for Estimates} \\
\hline
$k$ & \multicolumn{5}{|c|}{Scenario p-002}  & \multicolumn{5}{|c|}{Scenario p-016} \\
\hline
& $\mu_{k,0}$ & $\mu_{k,1}$ & $\sigma_{k,0}$ & $\sigma_{k,1}$ & $\xi_k$ & $\mu_{k,0}$ & $\mu_{k,1}$ & $\sigma_{k,0}$ & $\sigma_{k,1}$ & $\xi_k$ \\
\hline
1 & 25.58 & 0.67 & 0.27 & 0.00 & 0.13 & 17.38 & 0.39 & 0.00 & 0.00 & 0.11 \\
\hline
2 & 19.68 & 0.52 & 0.26 & 0.00 & 0.17 & 11.10 & 0.29 & 0.00 & 0.00 & 0.10 \\
\hline
3 & 14.89 & 0.43 & 0.28 & 0.00 & 0.18 & 7.42 & 0.23 & 0.00 & 0.00 & 0.09 \\
\hline
4 & 14.69 & 0.41 & 0.26 & 0.00 & 0.16 & 6.97 & 0.21 & 0.00 & 0.00 & 0.09  \\
\hline
5 & 13.42 & 0.37 & 0.29 & 0.00 & 0.20 & 7.60 & 0.19 & 0.00 & 0.00 & 0.12  \\
\hline
6 & 11.59 & 0.32 & 0.29 & 0.00 & 0.17 & 6.70 & 0.17 & 0.00 & 0.00 & 0.12 \\
\hline
7 & 8.84 & 0.26 & 0.32 & 0.00 & 0.13 & 4.82 & 0.15 & 0.00 & 0.00 & 0.15  \\
\hline
8 & 7.37 & 0.23 & 0.32 & 0.00 & 0.13 & 3.89 & 0.13 & 0.00 & 0.00 & 0.13 \\
\hline
9 & 7.91 & 0.22 & 0.32 & 0.00 & 0.12 & 4.27 & 0.12 & 0.00 & 0.00 & 0.14 \\
\hline
10 & 6.55 & 0.19 & 0.32 & 0.00 & 0.11 & 4.41 & 0.13 & 0.00 & 0.00 & 0.14 \\
\hline 
\end{tabular}
\caption{Maximum likelihood estimates for values of increasing window size using fixed shape for windows $k>1$.}\label{tab:nonstationaryGEV}
\end{table}

\begin{table}[h!]
\centering
\begin{tabular}{|c|c|c|}
\hline
\multicolumn{3}{|c|}{Extremal Index $\theta_k$ Approximation} \\
\hline
$k$ & \multicolumn{1}{|c|}{Scenario p-002}  & \multicolumn{1}{|c|}{Scenario p-016}\\
\hline
1 & 0.64 & 0.60 \\
\hline
2 & 0.60 & 0.59 \\
\hline
3 & 0.64 & 0.56 \\
\hline
4 & 0.60 & 0.51\\
\hline
5 & 0.60 & 0.48 \\
\hline
6 & 0.53 & 0.43 \\
\hline
7 & 0.53 & 0.42 \\
\hline
8 & 0.47 & 0.40 \\
\hline
9 & 0.43 & 0.38 \\
\hline
10 & 0.39 & 0.28 \\
\hline
\end{tabular}
\caption{Ferro-Segers estimate of the extremal index for increasing window size using nonstationary thresholding.}\label{tab:ei}
\end{table}

\begin{table}[h!]
\centering
\begin{tabular}{|c|c|c|c|}
\hline
\multicolumn{4}{|c|}{$g_T(k)$} \\
\hline
\multicolumn{2}{|c|}{Scenario p-002}  & \multicolumn{2}{|c|}{Scenario p-016}\\
\hline
$a$ & $b$ & $a$ & $b$ \\
\hline
0.80 & 0.80 & 0.86 & 0.79 \\
\hline
\end{tabular}
\caption{Approximations of the parameters in $g_T(k)$ for each scenario through generalised linear model fitting.}\label{tab:gkt}
\end{table}

\clearpage
\appendix
\section{Nonstationary fitting and related statistical tests}\label{append}

\begin{figure}[h!]
\centering
\includegraphics[width=\textwidth]{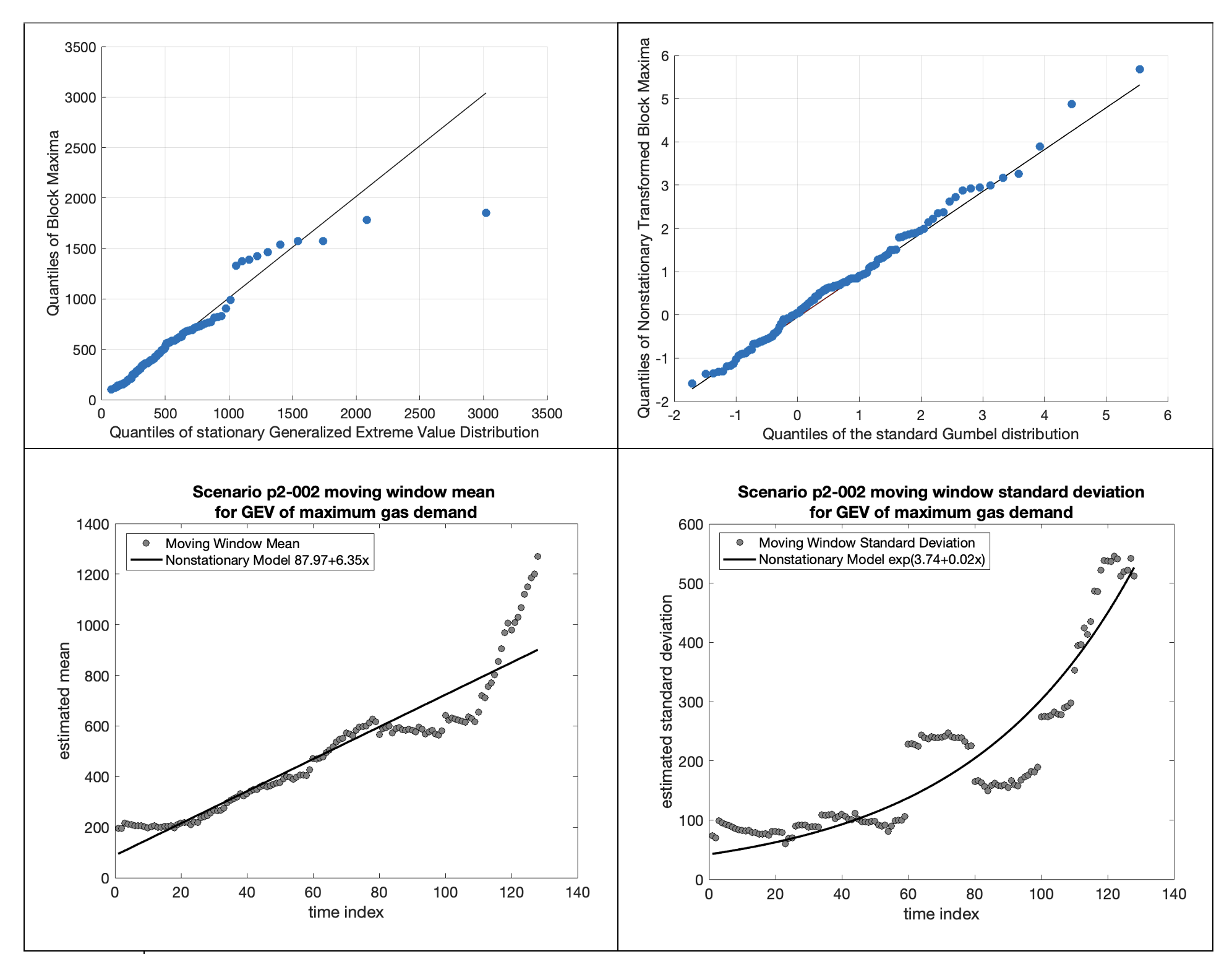}
\caption{(Top) Quantile plots before (left) and after (right) nonstationary fitting and (Bottom) empirically estimated mean and variance over time for scenario p2-002.}\label{fig:non}
\end{figure}

\begin{table}[h!]
\centering
\begin{tabular}{|c|c|c|c|c|}
\hline
\multicolumn{5}{|c|}{Statistical tests of fit across $k$ using standard MLE for each scenario} \\
\hline
 & \multicolumn{2}{|c|}{Scenario p-002}  & \multicolumn{2}{|c|}{Scenario p-016}\\
\hline
 $k$ & k-s test $p$-Value & a-d test $p$-Value & k-s test $p$-Value & a-d test $p$-Value\\
\hline
1 & 0.82 & 0.93 & 0.99 & 0.99 \\
\hline
2 & 0.85 & 0.89 & 0.96 & 0.99 \\
\hline
3 & 0.45 & 0.63 & 0.80 & 0.94  \\
\hline
4 & 0.96 & 0.97 & 0.86 & 0.91 \\
\hline
5 & 0.42 & 0.50 & 0.82 & 0.62 \\
\hline
6 & 0.53 & 0.66 & 0.50 & 0.42 \\
\hline
7 & 0.99 & 1.00 & 0.63 & 0.64 \\
\hline
8 & 0.86 & 0.92 & 0.79 & 0.90 \\
\hline
9 & 0.98 & 0.98 & 0.98 & 0.99 \\
\hline
10 & 0.77 & 0.84 & 0.55 & 0.90 \\
\hline
\end{tabular}
\caption{$p$-Values for the Kolmogorov-Smirnov (k-s) test and Anderson-Darling (a-d) test of fit on the data using standard maximum likelihood estimation of the derived non-stationary model. The resulting approximated parameters are used to estimate the function $g_T(k)$ which is used in the proposed model to estimate the distributions for $k$ beyond $k_f$, the successive window for which the MLE fails.}\label{tab:fit}
\end{table}

\clearpage
\bibliography{ref}

\end{document}